\documentclass[12pt]{article}
\usepackage{a4}
\usepackage{latexsym}
\usepackage{amsmath}
\usepackage{amssymb}
\usepackage{amscd}
\usepackage{paralist}
\usepackage{color}
\usepackage{amsthm}
\theoremstyle{plain}
    \newtheorem{theorem}                    {Theorem}       [section]
    
    \newtheorem{corollary}  [theorem]       {Corollary}
    \newtheorem{proposition}[theorem]       {Proposition}

\newtheorem{remark}[theorem]{Remark}

\DeclareFontFamily{U}{wncy}{}
    \DeclareFontShape{U}{wncy}{m}{n}{<->wncyr10}{}
    \DeclareSymbolFont{mcy}{U}{wncy}{m}{n}
    \DeclareMathSymbol{\Sh}{\mathord}{mcy}{"58}

\newcommand{\Hom}{\operatorname{Hom}}

\newcommand{\Pic}{\operatorname{Pic}}

\newcommand{\Spec}{\operatorname{Spec}}

\newcommand{\im}{\operatorname{im}}

\newcommand{\coker}{\operatorname{coker}}

\newcommand{\Br}{\operatorname{Br}}

\renewcommand{\lim}{\operatorname{lim}}

\newcommand{\f}{{\mathcal F}}
\newcommand{\g}{{\mathcal G}}

\newcommand{\Z}{{{\mathbb Z}}}

\newcommand{\Q}{{{\mathbb Q}}}
\newcommand{\R}{{{\mathbb R}}}
\newcommand{\F}{{{\mathbb F}}}
\newcommand{\G}{{{\mathbb G}_m}}

\newcommand{\X}{{X}}

\newcommand{\et}{{\text{\rm et}}}

\newcommand{\sep}{{\text{\rm s}}}

\renewcommand{\O}{{\cal O}}

\newcommand{\proofend}{\hfill$\square$\\ \smallskip}

\date{}

\title{Comparing the Brauer group to the Tate-Shafarevich group.}

\author{Thomas H. Geisser}

\begin{document}

\renewcommand{\thefootnote}{\fnsymbol{footnote}} 
\footnotetext{\emph{MSC classes}: 11G40, 14G17, 14J20, 11G25, 11G35}

\maketitle

\begin{abstract}
We give a formula relating the order of the Brauer group of a surface
fibered over a curve over a finite field to the order of the
Tate-Shafarevich group of the Jacobian of the generic fiber. 
The formula implies that the Brauer group of a smooth and proper 
surface over a finite field is a square if it is finite.
\end{abstract}

\section{Introduction}
Let $K$ be a global field, and let 
$V$ be the smooth and proper model if $K$ has characteristic $p$,
or the spectrum of the ring of integers of $K$ in the number field case.
Let $\X$ be a regular surface and $\X\to V$ a projective flat map 
with geometrically connected fibers such that $X_K=\X\times_VK$ 
is smooth over $K$. 
For point $v\in V$, let $K_v$ be the completion of $K$
and $X_{K_v}=X\times_VK_v$. 

It is a classical result of Artin
and Grothendieck \cite{brauerIII} that the Brauer group of $\X$ is finite 
if and only if the Tate-Shafarevich group of the Jacobian $A=\Pic^0_{X_K}$ of $X_K$ 
is finite.
Grothendieck \cite[(4.7)]{brauerIII}, Milne \cite{milne}, 
and Gonzalez-Aviles \cite{gonzalez} gave formulas relating the order
of the Brauer group of $\X$ to the order of the Tate-Shafarevich group $\Sh(A)$
under some conditions on the periods of $X_{K_v}$. We give a 
general formula without any conditions. 
Let $\delta$ and $\delta_v$ be the index of $X_K$ and $X_{K_v}$,
respectively, and $\alpha$ and $\alpha_v$ be the order of 
the cokernel of the inclusion $\Pic^0(X_K)\to H^0(K,\Pic_{X_K}^0)$
and $\Pic^0(X_{K_v})\to H^0(K_v,\Pic_{X_K}^0)$, respectively. 
By Lichtenbaum \cite[Thm.3 (proof)]{lichtenbaumcurve}, 
$\alpha_v$ is equal to the period 
$\delta_v'$ of $X_{K_v}$. 

\begin{theorem}\label{mmm}
If $K$ has no real embeddings and if the Brauer group $\Br(\X)$ is finite, then 
\begin{equation}\label{lformul}
|\Br(\X)|\alpha^2  \delta^2= |\Sh(A)|\prod_{v\in V} \alpha_v\delta_v.
\end{equation}
\end{theorem}

This generalizes the results of Grothendieck, Milne and Gonzalez-Aviles, 
and corrects the formula of Liu, Lorenzini and Raynaud \cite{liu}
by the factor $\alpha^2$. The problem is that \cite{liu} uses the incorrect
\cite[Lemma 4.2]{gordon}, which implies that $\alpha=1$, see their corrigendum
\cite{liucor}.
If $K$ is a number field with real embeddings, then the same formula holds
up to a power of $2$ (due to the usual problem with duality for
Galois cohomology of a number ring with real places).
By \cite[Remark 4.5]{llr}, the right hand side in \eqref{mmm}
is a square, hence the argument of \cite{liu}
gives the following

\begin{corollary}
Let $X$ be a smooth and proper surface
over a finite field. If the Brauer group is finite, then its order
is a square.
\end{corollary}

A key ingredient in the proof is the following local-to-global result for
the Brauer group:

\begin{theorem}\label{cokerl}
If $\Br(\X)$ is finite and $K$ has no real embeddings, then
$$ 0\to \Br(\X)\to \Br(X_K)\to \bigoplus_{v\in V} \Br(X_{K_v})
\to \Hom(\Pic(X_K),\Q/\Z) \to 0$$
is exact.
\end{theorem}


\smallskip 

We thank T. Szamuely and C. Gonzalez-Aviles
for comments on an earlier version of this paper, and the referee
for his careful reading.

\section{Brauer groups and Tate-Shafarevich groups}
We continue to use the notation of the introduction. 
For a closed point $v$ of $V$, we let $\O_v$ be the completion of 
$V$ at $v$,  $k_v$ the residue field at $v$, 
and $Y_v=\X\times_V k_v$.
Let $G$ and $G_v$ be the Galois group of $K$ and $K_v$, respectively.

Denoting the Pontrjagin dual of the abelian group $A$ by 
$A^*=\Hom(A,\Q/\Z)$, we have Lichtenbaum's duality for the curve $X_{K_v}$ 
\cite{lichtenbaumcurve}
\begin{equation}\label{lpair}
\Pic(X_{K_v})^* \cong \Br(X_{K_v}).
\end{equation}
This duality has been generalized by S.\ Saito to include the finite
characteristic case in \cite[Thm. 9.2]{saitod}. Both Lichtenbaum's
and Saito's pairing are defined by pulling back elements of 
$\Br(X_{K_v})$ along divisors, and checking that the result vanishes
on principal divisors.
Composing with the dual of the natural map
$\Pic(X_K)\to \prod_{v\in V}\Pic(X_{K_v})$,
we obtain a map of discrete torsion groups
$$ \bigoplus_{v\in V} \Br(X_{K_v})\stackrel{l}{\to} \Pic(X_K)^*.$$


\smallskip

{\it Proof of Theorem \ref{cokerl}:}
Exactness at the left two terms 
can be found in \cite[Lemma 2.6]{milne}. Exactness
on the right follows from \eqref{lpair} and injectivity of
$\Pic(X_K)\to \Pic(X_{K_v})$ for every $v$. Indeed, if $X_K$
has a point over a finite Galois extension $L$ and $w$ is a place of $L$
above $w$, then $\Pic(X_K)\to \Pic(X_L)\to \Pic(X_{L_w})$ is injective,
the former by the Hochschild-Serre spectral sequence and the latter because
the Picard functor is representable in the presence of a rational point.
It remains to show exactness at the sum. Consider the diagram
$$\begin{CD}\label{mmdd}
\Br (X_K) @>>>  \bigoplus_{v\in V} \Br(X_{K_v})
@>>> \Pic(X_K)^* @>inj >> \Pic(\X)^* \\
@| @Vinj VV @. @|\\
\Br (X_K) @>>> \bigoplus_{v\in V} H^{3}_{Y_v}(\X,\G)
@>>> H^3_\et(\X,\G) @>\xi>> \Pic(\X)^*.
\end{CD}$$ 
The left three terms of the second row 
arise from the localization sequence for etale cohomology for $\X$, and the
second vertical injection is the sum of the localization sequences 
for the $\X_{\O_v}$, using the vanishing of $\Br(\X_{\O_v})$ 
(\cite[Thm. 3.1]{brauerIII}, see the
proof of \cite[Lemma 2.6]{milne}), and the fact that 
$H^{3}_{Y_v}(\X,\G) \cong  H^{3}_{Y_v}(\X_{\O_v},\G)$. 
A diagram chase shows that the exactness at the sum follows if we can 
define an injective map $\xi$ such that the right rectangle commutes.

We define the map $\xi$ by using a divisor $D$ on $X$ 
to pull-back cohomology classes in $H^3_\et(\X,\G)$ to the normalization 
$H^3_\et(\tilde D,\G)$, which is isomorphic to $(\Q/\Z)^c$, 
$c$ the number of irreducible components of $D$
\cite[II Rem. 2.2 (b)]{adt}, and then summing up. Then the right rectangle commutes because both
composition are defined by pulling back cohomology classes along divisors.
Saito defines a map $\phi^1: H^3_\et(\X,\G) \to  \Pic(\X)^*$ and shows in
\cite[Thm. 5.5(2)]{saito} that it is a surjection whose kernel 
vanishes if $\Br(\X)$ is finite. It suffices to show that $\phi^1=\xi$. 
In the proof of loc.cit., one chooses a divisor $Y$ whose components generate $\Pic(X)$, 
pulls-back cohomology classes $H^3_\et(X,\G)$ to $H^3_\et(Y,i^*\G)$, and uses the duality 
between $H^3_\et(Y,i^*\G)\cong (\Q/\Z)^c$ and
$H^1_Y(X,\G)\cong \Z^c$, where $c$ is the number of components of $Y$
\cite[Prop. 4.6]{saito}. Now it suffices to observe that under the given hypothesis,
the map $H^3_\et(X,\G)$ to $H^3_\et(Y,i^*\G)$ is injective, 
the map $\Z^c \cong H^1_Y(X,\G)\to \Pic(X)$ sends a generator corresponding
to a component of $Y$ to its divisor class, and 
$$H^3_\et(Y,i^*\G)\cong  H^3_\et(Y,\G)\cong H^3_\et(\tilde Y,\G)
\cong (\Q/\Z)^c,$$
which follows from the proof of \cite[(4-11)]{saito}.
\proofend

\begin{remark}
1) In the function field case one can show that the sequence
is exact except at the sum, where its cohomology is $(T\Br(\X))^*$
up to $p$-groups. 

2) The hypothesis on $2$-torsion in case of real embeddings is used
to apply S.Saito's result, see \cite[\S 5]{saito}.
\end{remark}

The following generalization of the Cassels-Tate exact
sequence by Gonzalez-Aviles and Tan \cite{tan}
can be thought of as the analog for the Tate-Shafarevich group.
The results on flat cohomology that are used have been corrected
in \cite{harari}.

\begin{theorem}\label{rightcol}
Let $A$ be an abelian variety over $K$ with dual $A^t$, and assume that
the Tate-Shafarevich group $\Sh(A^t)$ is finite. Then the sequence
$$0\to  \Sh(A)\to H^1(K,A )\stackrel{\beta^1}{\to} \bigoplus_v H^1(K_v,A)
\stackrel{\gamma^1}{\to} H^0(K,A^t)^*\to 0$$
is exact. 
\end{theorem}

Here the map $\gamma^1$ is the dual of the injection
$$\beta^0: H^0(K,A^t)^\wedge \to 
\prod_v H^0(K_v,A^t)^\wedge\cong (\bigoplus_v H^1(K_v,A))^*,$$
where $G^\wedge=\lim_m G/m$ denotes the completion of an abelian group $G$.

\section{Comparison}
We complete the proof of Theorem \ref{mmm} by comparing the sequences of Theorem \ref{cokerl}
and of Theorem \ref{rightcol} applied to $\Pic^0_X$ via their maps to 
$$H^1(K,\Pic_{X_K})\to \bigoplus_v H^1(K_v,\Pic_{X_K}).$$

The long exact sequence of Galois cohomology groups
associated to the degree map over the separable closure $K^s$ of $K$
$$0\to \Pic(X_{K^\sep})^0\to 
\Pic(X_{K^\sep})\stackrel{\deg}{\longrightarrow} \Z\to 0$$
induces the middle two exact rows of the following diagram:
\begin{equation}\label{dia2}
\begin{CD}
@.  @. \Sh(\Pic_{X_K}^0)@>>> \Phi \\
@. @. @VVV @VVV \\
0@>>> \Z/\delta' @>>> H^1(K,\Pic_{X_K}^0 )@> >>  H^1(K,\Pic_{X_K})@>>> 0 \\
@.@VVV @V\beta^1 VV @V\tau VV \\
0@>>> \oplus \Z/\delta_v'@>>> \bigoplus_v H^1(K_v,\Pic^0_{X_K} )@>>> 
\bigoplus_v H^1(K_v,\Pic_{X_K})@>>> 0\\
@. @.  @V\gamma^1 VV @V\rho VV \\
@. @.  H^0(K,\Pic^0_{X_K})^*@>\omega>> \Psi \\
 \end{CD}\end{equation}
The upper and lower rows are the kernels and cokernels of the vertical maps.
Finiteness of $\Sh(\Pic_{X_K}^0)$ and of $\oplus \Z/\delta_v'$ implies finiteness of $\Phi$. 
Counting orders we obtain the formula
$$ |\Phi| = 
\frac{ |\Sh(\Pic_{X_K}^0)|\cdot \prod_v\delta_v'}{|\ker \omega|\cdot \delta'}.$$

Now we use the (functorial) Hochschild-Serre spectral sequence 
\begin{equation}\label{fhss}
0\to \Pic(X_K)\to H^0(K,\Pic_{X_K})\to   \Br(K)\to   \Br(X_K) \to
H^1(K,\Pic_{X_K})\to 0
\end{equation}
for $X$ and $X_{K_v}$ 
to obtain the middle two exact rows of the following diagram:
\begin{equation}\label{eins+}\begin{CD} 
@. 0@>>> \Br(X) @>>> \Phi\\
@.@VVV @VVV @VVV\\
0\to\;\;  P\qquad @>>> \Br(K) @>>>\Br(X_K) @>>> 
H^1(K,\Pic_{X_K}) \quad\; \to 0\\
@VVV @V VV @VVV @V\tau VV \\
0\to  \bigoplus \Z/\delta_v@>>> \bigoplus \Br(K_v)@>>>\bigoplus  \Br(X_{K_v})@>>> 
\bigoplus H^1(K_v,\Pic_{X_K})\to 0\\
@. @V\sum VV @VVV@V\rho VV \\
 @. \Q/\Z @>\deg^*>>  \Pic(X_K)^* @>\sigma>>\qquad \qquad\;\Psi\qquad \; \to 0
\end{CD}
\end{equation} 
The upper and lower rows are the kernels and cokernels of the vertical maps. 
The kernel of $\Br(K_v)\to \Br(X_{K_v})$ is isomorphic to  $\Z/\delta_v$
by the Lichtenbaum-Roquette theorem \cite[Thm. 3]{lichtenbaumcurve}.
Lichtenbaum's result is stated in characteristic $0$,
but the proof works in characteristic $p$ as soon as 
duality for Galois cohomology of abelian varieties holds 
\cite[Thms. 9.2, 9.3]{saitod}.  
The lower middle square is commutative by the definition of 
the pairing \eqref{lpair},
see \cite[p. 125]{lichtenbaumcurve}, and functoriality of the degree map:
$$\begin{CD}
\Br(K_v)@= \Q/\Z @= \Q/\Z\\
@VVV @V\deg^*VV @V\deg^* VV\\
\Br(X_{K_v})@>>> \Pic(X_{K_v})^* @>>> \Pic(X_K)^*.
\end{CD}$$
We have $|\ker \deg^*|=\delta$, and 
since $\coker \deg^*\cong \Pic^0(X_K)^*$, $\sigma$ factors through a map 
$\zeta :  \Pic^0(X_K)^*\to \Psi$, and $\ker \sigma/\im deg^*\cong \ker \zeta$.
Counting orders we obtain the formula
$$ |\Phi|= \frac{ |\Br(X)|\cdot |P|\cdot\delta}
{|\ker \zeta| \cdot \prod_v\delta_v}.$$
In order to relate the maps $\omega$ and $\zeta$, we consider the following 
diagram. 
$$\begin{CD}
\Br(X_{K_v}) @>surj>> H^1(K_v,\Pic_{X_K})@<surj << H^1(K_v,\Pic^0_{X_K})\\
@|@|@| \\
\Pic(X_{K_v})^* @>surj>> \Pic^0(X_{K_v})^*@<surj << H^0(K_v,\Pic^0_{X_K})^*\\
@VVV @VVV @V\gamma^1VV \\
\Pic(X_K)^*@>surj>> \Pic^0(X_K)^*@<f^*<surj< H^0(K,\Pic^0_{X_K})^*\\
@| @V\zeta VV@| \\
\Pic(X_K)^*@>\sigma >> \Psi @<\omega<< H^0(K,\Pic^0_{X_K})^*
\end{CD}$$
If we replace the middle composition by $\rho$, then it defines 
the maps $\sigma$ of \eqref{eins+} in the left half and $\omega $ 
of \eqref{dia2} in the right half.
Since $\sigma$ factors
through $\zeta$, so does $\rho$. 
The upper two squares are commutative by compatibility
of Lichtenbaum's perfect pairings \cite[\S 4]{lichtenbaumcurve},
and the other squares are obviously commutative. Then the left half
of the diagram shows that the middle composition as indicated agrees
with the the map $\rho$, and the right half of the diagram shows that 
$\omega=\zeta f^*$. Since $f^*$ is surjective, we obtain
$$ |\ker \omega|=|\ker f^*|\cdot |\ker \zeta|.$$
Finally, the diagram 
\begin{equation}\label{diagr}\begin{CD}
0@>>> \Pic^0(X_K) @>>> \Pic(X_K)@>>> \delta\Z@>>> 0\\
@.@VfVV @VVV @VVV \\
0@>>> H^0(K,\Pic^0_{X_K}) @>>> H^0(K, \Pic_{X_K})@>>> \delta'\Z@>>> 0
\end{CD}
\end{equation}
shows that $|\ker f^*|\cdot \delta = |P|\cdot \delta'$. Since 
$\alpha=|\ker f^*|$, we obtain Theorem \ref{mmm} by equating the two formulas for $|\Phi|$.

\begin{remark}
The following example that for a curve $C$ over a global field
$K$, the $l$-rank of $\coker \Pic^0(C)\to H^0(K,\Pic^0_C)$ 
can be arbitrary large for any $l$ was communicated
to us by J.\ Starr. By the sequence \eqref{fhss} and the diagram 
\eqref{diagr}, it suffices to find $C$ such that the $l$-rank of  
the kernel of $\Br(K)\to \Br(C)$ is arbitrary large. 
Let $a_0, \ldots ,a_r$ be $\Z/l$-linearly independent classes 
in $Br(K)[l]$, and let $P_0,\ldots,P_r$ be the associated 
Severi-Brauer $K$-schemes. Now let $C$ be a general complete 
intersection curve in the product variety
$P=P_0\times_K\cdots \times_KP_r$. Then
the kernel of the pullback map $\Br(K)\to \Br(P)$ contains the classes
$a_0,\ldots, a_r$, hence so does the kernel of $\Br(K)\to \Br(C)$.
\end{remark}

\small

{\sc  Rikkyo University, Ikebukuro, Tokyo, Japan}

\textit{E-mail address:} {\tt geisser@rikkyo.ac.jp}
\end{document}